%% file: ifacconf.tex
\documentclass{ifacconf}

\usepackage{amsmath,amsfonts}

\usepackage{graphicx}      

\makeatletter
\let\old@ssect\@ssect 
\makeatother

\usepackage{natbib}        

\usepackage{scalerel}
\usepackage{tikz}
\usetikzlibrary{svg.path}
\input{ORCID_logo}

\usepackage{hyperref}

\makeatletter
\def\@ssect#1#2#3#4#5#6{%
	\NR@gettitle{#6}
	\old@ssect{#1}{#2}{#3}{#4}{#5}{#6}
}
\makeatother

\begin{document}
	
\begin{frontmatter}

\title{Bifurcation preserving discretisations of optimal control problems} 

\thanks[copyright]{\textcopyright 2021 the authors. This work has been accepted to IFAC for publication under a Creative Commons Licence CC-BY-NC-ND}

\author[First]{Christian Offen}
\author{and Sina Ober-Blöbaum$^\ast$}

\address[First]{Department of Mathematics, Paderborn University, Germany\\ 
(e-mail: christian.offen@uni-paderborn.de) }

\begin{abstract}                
	The first order optimality conditions of optimal control problems (OCPs) can be regarded as boundary value problems for Hamiltonian systems. Variational or symplectic discretisation methods are classically known for their excellent long term behaviour. As boundary value problems are posed on intervals of fixed, moderate length, it is not immediately clear whether methods can profit from structure preservation in this context. When parameters are present, solutions can undergo bifurcations, for instance, two solutions can merge and annihilate one another as parameters are varied. We will show that generic bifurcations of an OCP are preserved under discretisation when the OCP is either directly discretised to a discrete OCP (direct method) or translated into a Hamiltonian boundary value problem using first order necessary conditions of optimality which is then solved using a symplectic integrator (indirect method). Moreover, certain bifurcations break when a non-symplectic scheme is used. The general phenomenon is illustrated on the example of a cut locus of an ellipsoid.
\end{abstract}


\end{frontmatter}

\input{CustomDefs}

\section{Introduction}\label{sec:intro}

%
There are two main strategies to discretise optimal control problems (OCPs): {\em direct} and {\em indirect} methods. 
In direct methods the OCP is approximated by a discrete optimisation problem, which is then solved using techniques from nonlinear programming. In indirect methods first order necessary conditions for optimality are calculated for the OCP. These have the structure of a boundary value problem for a Hamiltonian system (Pontryagin's principle). The boundary value problem is solved numerically using methods such as shooting algorithms or implicit solvers for the fully discretised problem. 
For this, Hamilton's equations need to be discretised. If a symplectic partitioned Runge-Kutta method is used, then the scheme is mathematically equivalent to a direct method, in which the state equation is integrated with the underlying Runge-Kutta method. If, on the other hand, a non-symplectic integrator is used, then the scheme cannot be obtained as a direct discretisation method \citep{ober2011DMOC}.

While there is some mathematical beauty in the fact that forming first order necessary conditions for optimality and discretisation commutes, provided that a symplectic integration scheme is used, any practical relevance of structure preservation in this context may not be immediately clear. Indeed, it has been argued that for OCPs symplectic integrators have no advantages over non-symplectic schemes with the exception of some special cases \citep{SymplecticityOptimalControl}. However, in this paper we show that using structure preserving integration schemes can be crucial when {\em bifurcation phenomena} occur.


Solutions to first order necessary conditions of parameter-dependent OCPs may not be unique but bifurcate as parameters are varied: for instance, two solutions can merge and annihilate one another or three solutions can interact.
More specifically, we will focus on families of OCPs where a cost function
\begin{equation}
	\label{eq:OCPGeneralIntro}
	S(u;\mu)= \int_{t_0}^{t_N}L(q(t),u(t);\mu) \d t
\end{equation}
is extremised subject to a state equation
\begin{equation}
	\label{eq:stateEQ}
\dot{q} = f(q,u;\mu), \quad q(t_0)=q_0(\mu), \, q(t_N)=q_N(\mu)
\end{equation}
among all admissible controls $u$. Here $\mu \in \Lambda$ is the parameter of the family of OCPs and $\Lambda$ is the parameter space. The parameter is fixed during optimisation.


An analysis of the bifurcation behaviour of solutions to OCPs helps to determine for which parameter values a unique optimal solution exists and in which parameter ranges there are several solutions which fulfil first order optimality conditions.
These bifurcation phenomena should be contrasted to bifurcations analysed in the literature related to branching due to low regularity \citep{Kogan1986}. Here, we restrict to a description of bifurcation phenomena of regular solutions which do not interact with boundaries of the state or control space.

\begin{figure}
	\begin{center}
		\includegraphics[width=0.3\linewidth]{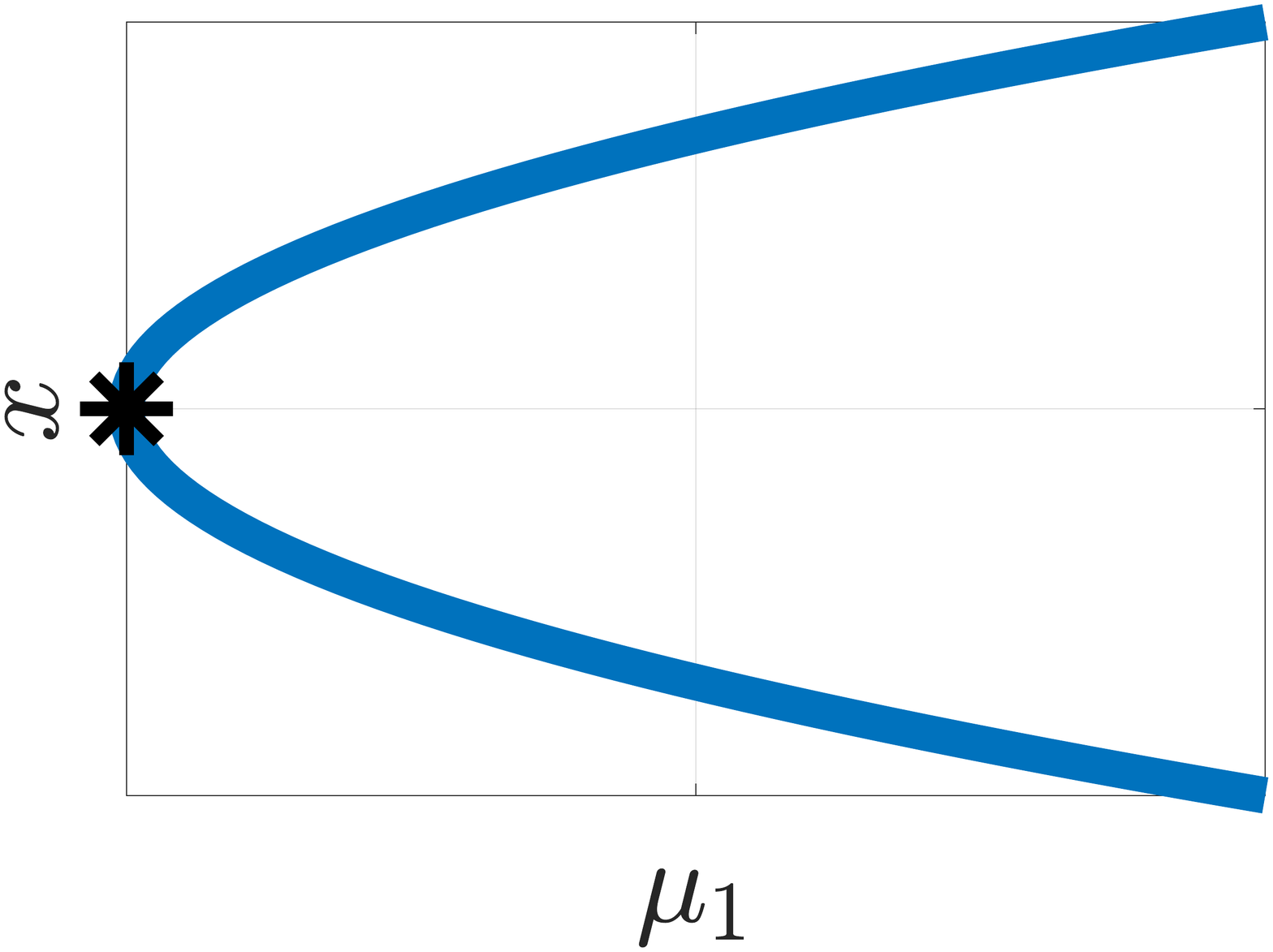}
		\includegraphics[width=0.3\linewidth]{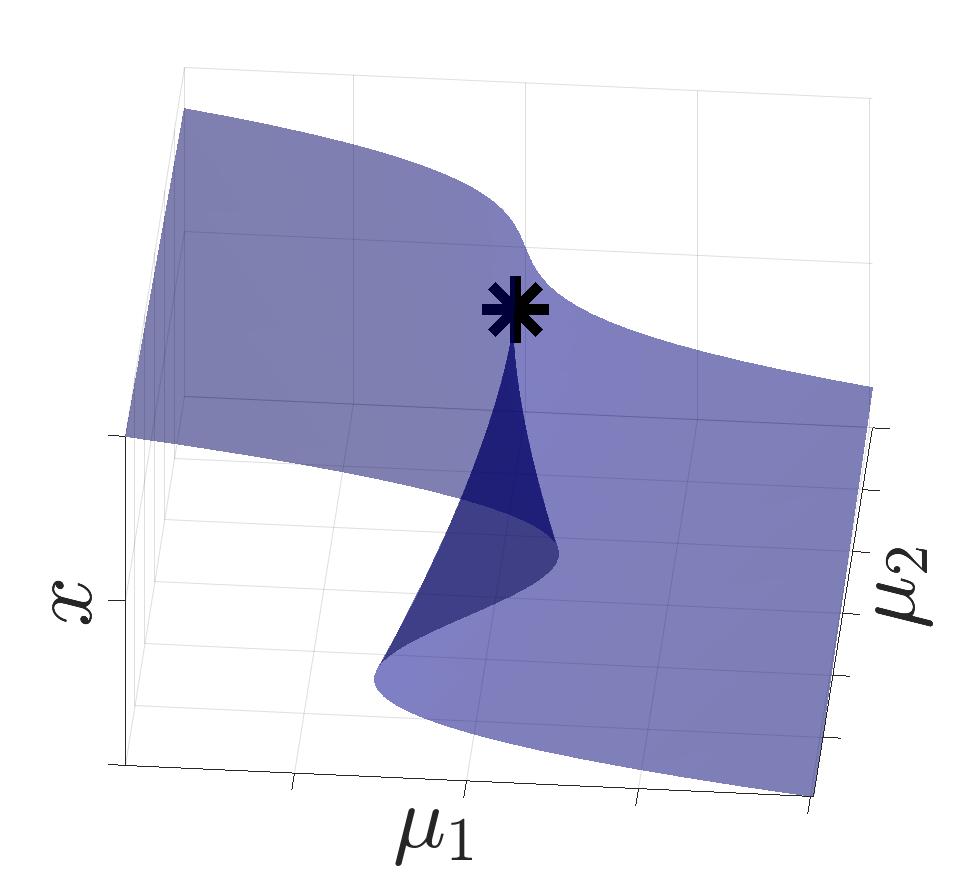}
		\includegraphics[width=0.3\linewidth]{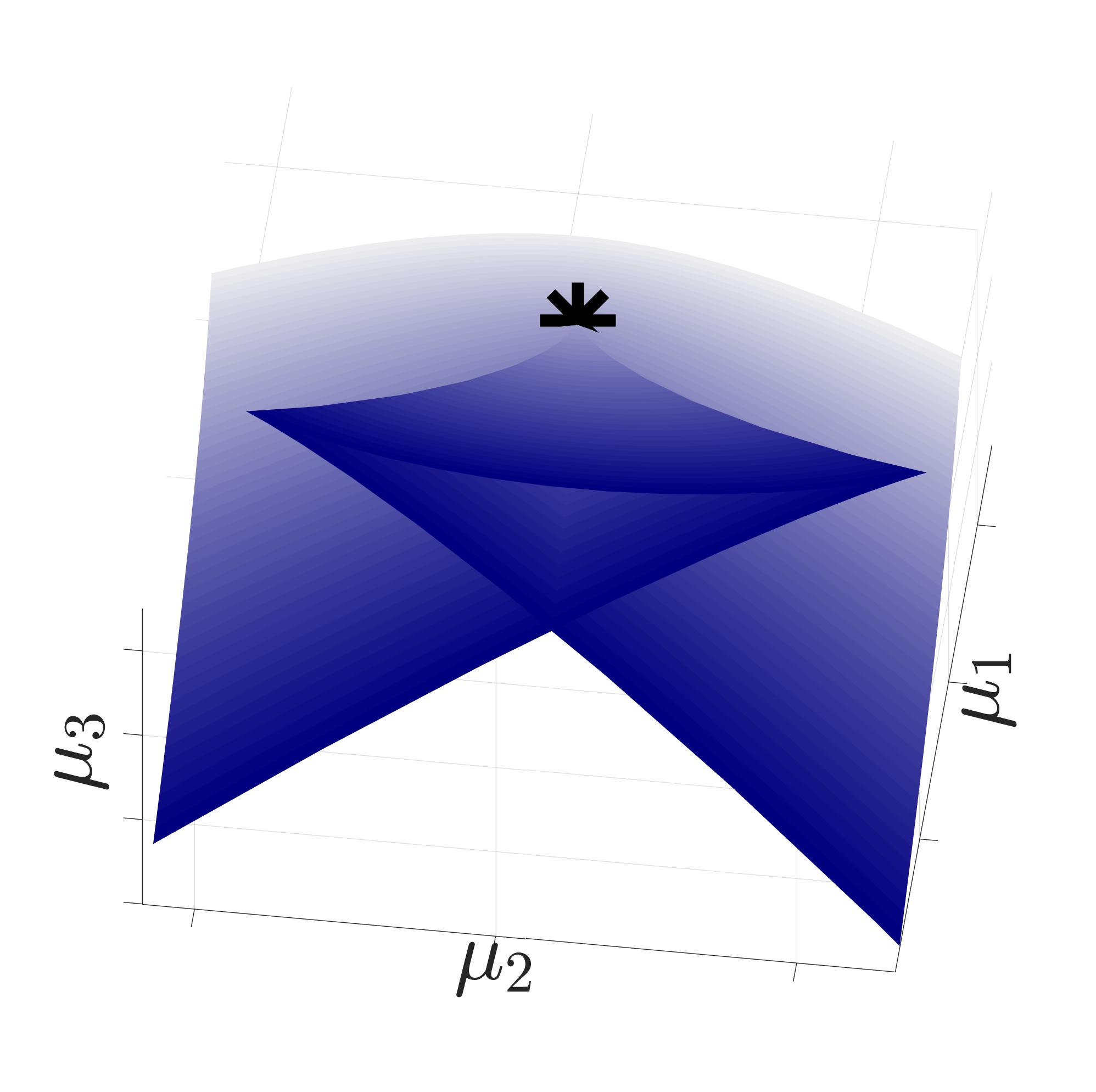}
		\includegraphics[width=0.3\linewidth]{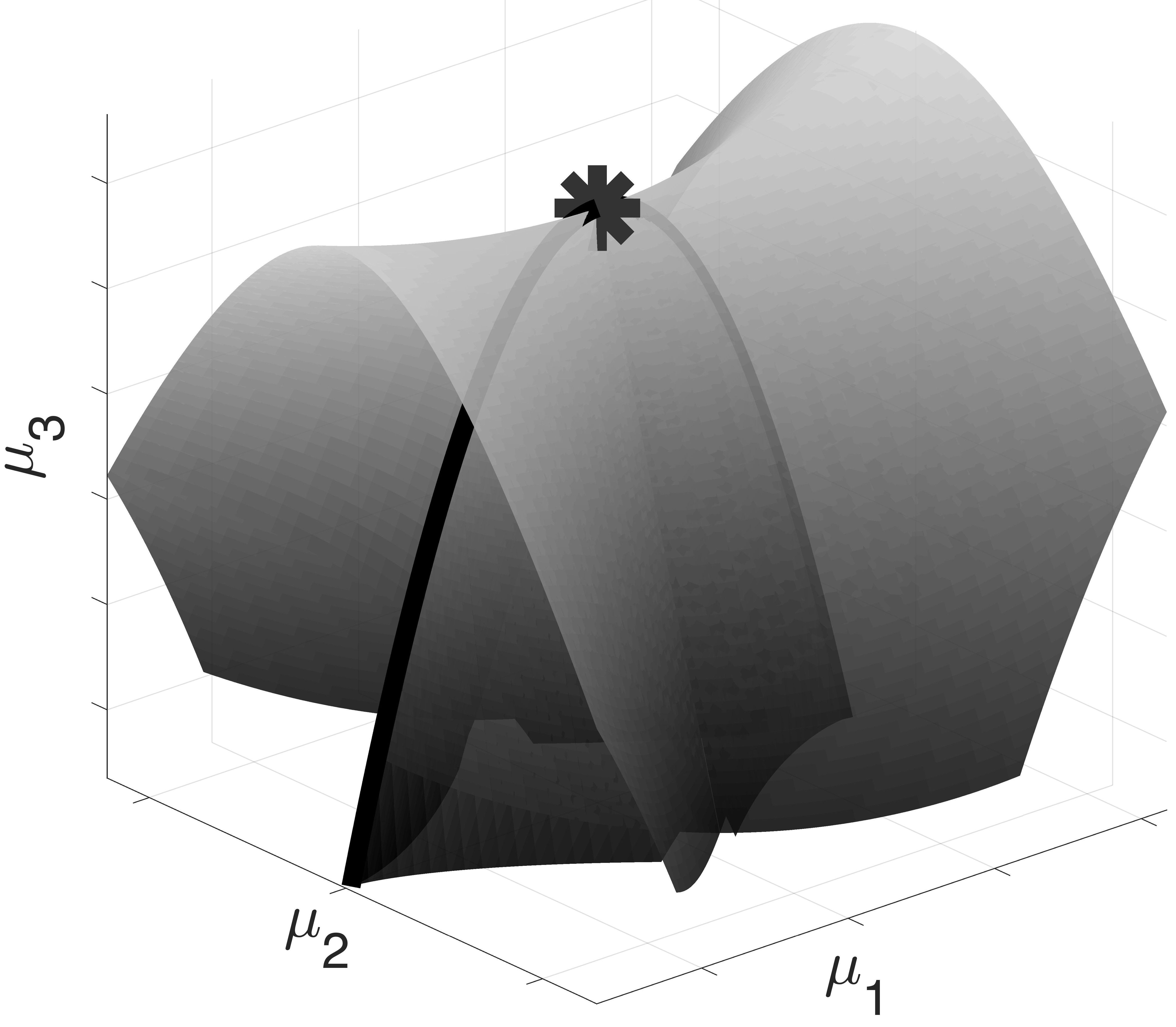}
		\includegraphics[width=0.3\linewidth]{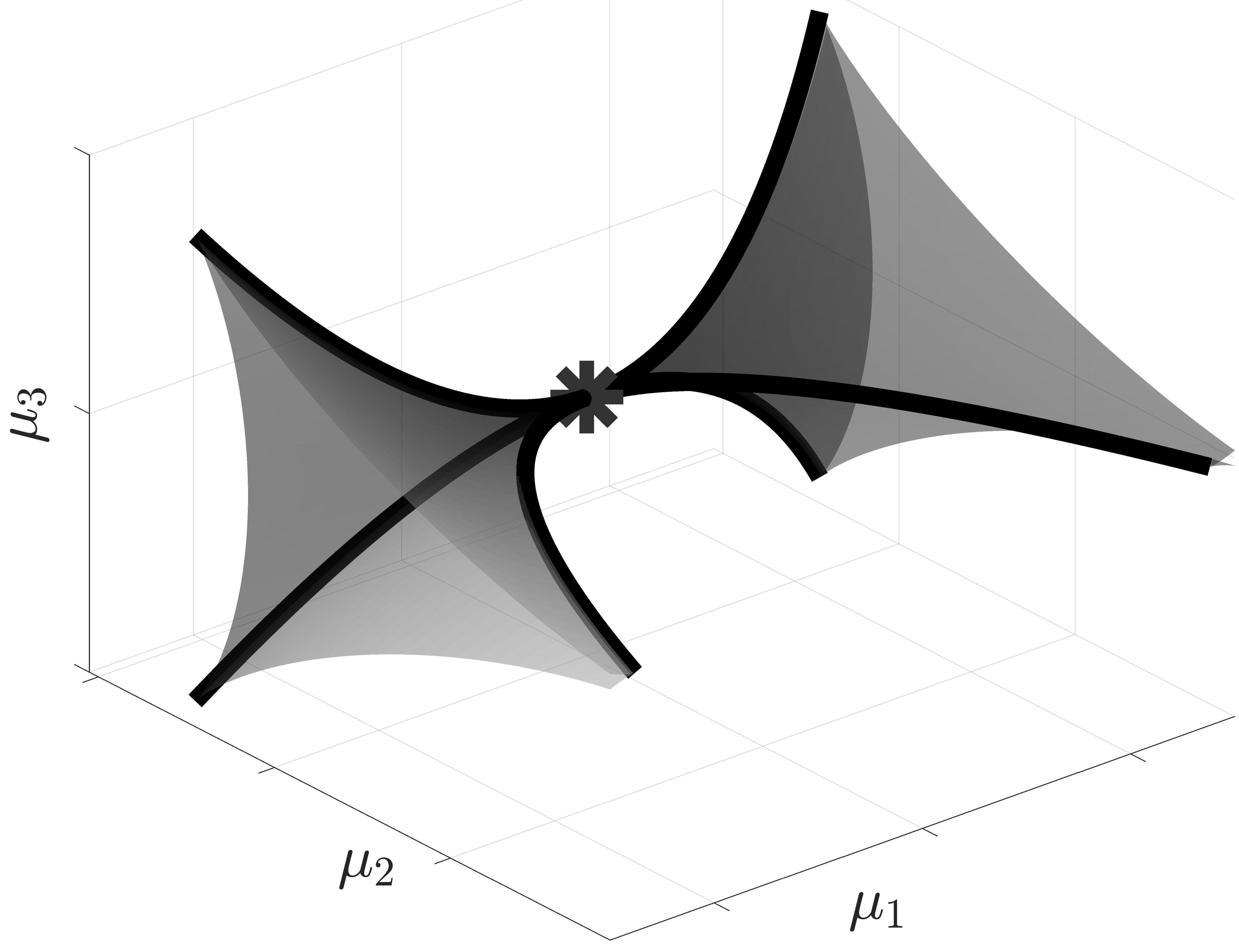}
	\end{center}
	\caption{The plots show generic behaviour of critical points of functions $r_\mu \colon \R^k \to \R$ when parameters $\mu$ are varied. From left to right, top to bottom we see models for fold, cusp, swallowtail, hyperbolic umbilic, and elliptic umbilic bifurcations. The most singular points are marked by $\ast$. They are persistent under small perturbations of $r_\mu$. See \cite{YoutubeSingularitiesAnimations} for animations.}\label{fig:Cat}
\end{figure}

The main example of the paper is the computation of shortest paths connecting two points $q_0$ and $q_N$ on an ellipsoid. There can be several connecting curves which extremise the length functional on an ellipsoid. These are geodesics. We will formulate the problem of finding geodesics starting at $q_0$ as a family of optimal control problems, where the parameter $\mu$ is given by the end point $q_N$. As $q_N$ is varied while $q_0$ is kept fix, the solutions bifurcate. The set of $q_N$, for which a bifurcation takes place, forms the {\em cut locus to $q_0$}. We will use popular discretisation methods from optimal control theory to compute cut loci and demonstrate that direct methods as well as indirect methods in combination with symplectic integrators resolve the loci correctly while non-symplectic integrators make qualitative errors.

Classical catastrophe theory considers the bifurcation behaviour of critical points of smooth, parameter dependent functions $r_\mu \colon \R^k \to \R$, where $\mu \in \R^l$ is a parameter. Stable bifurcations have been classified, see, for instance \citep{Arnold1}. Figure \ref{fig:Cat} shows the first five elementary catastrophes.
We relate bifurcations of solutions to first order necessary conditions of OCPs to classical catastrophe theory and explain why only direct methods and indirect methods in combination with symplectic integrators can preserve all stable bifurcations under discretisation.

The article is structured as follows. To exemplify the theoretical aspects of this work and to show their significance for optimal control theory, we introduce the main example of the paper, the computation of cut loci, and provide an optimal control formulation of the problem in section \ref{sec:GeoOC}. In section \ref{sec:Disc} we apply typical discretisation schemes to the optimal control formulation. Section \ref{sec:Ex} contains a numerical experiment, in which the conjugate locus of an ellipsoid is computed with the different discretisation schemes. The experiment demonstrates that some discretisation approaches preserve important qualitative aspects of the conjugate locus, while others break them. That the observed behaviour is prototypical for a large class of OCPs and discretisation schemes is proved in section \ref{sec:Theo}. The section, furthermore, connects bifurcations in OCPs to one of the authors' classification results for bifurcations in Hamiltonian boundary value problems \citep{bifurHampaper,numericalPaperShort,numericalPaper,PhDThesis}, on which this work is based. Section \ref{sec:Conclusion} summarises the findings.

\input{sec_geodesics}

\input{sec_experiment}

\input{sec_general}



\vfill

{ ORCID}

\noindent
Christian Offen:	\quad \quad 0000-0002-5940-8057 {\protect \orcidicon{0000-0002-5940-8057}}\\
Sina Ober-Blöbaum:  0000-0001-6720-7493 {\protect \orcidicon{0000-0001-6720-7493}}

\pagebreak

\bibliography{bibliography_doi.bib}             

\end{document}

%% file: ORCID_logo.tex
\definecolor{orcidlogocol}{HTML}{A6CE39}
\tikzset{
  orcidlogo/.pic={
    \fill[orcidlogocol] svg{M256,128c0,70.7-57.3,128-128,128C57.3,256,0,198.7,0,128C0,57.3,57.3,0,128,0C198.7,0,256,57.3,256,128z};
    \fill[white] svg{M86.3,186.2H70.9V79.1h15.4v48.4V186.2z}
                 svg{M108.9,79.1h41.6c39.6,0,57,28.3,57,53.6c0,27.5-21.5,53.6-56.8,53.6h-41.8V79.1z M124.3,172.4h24.5c34.9,0,42.9-26.5,42.9-39.7c0-21.5-13.7-39.7-43.7-39.7h-23.7V172.4z}
                 svg{M88.7,56.8c0,5.5-4.5,10.1-10.1,10.1c-5.6,0-10.1-4.6-10.1-10.1c0-5.6,4.5-10.1,10.1-10.1C84.2,46.7,88.7,51.3,88.7,56.8z};
  }
}

\newcommand\orcidicon[1]{\href{https://orcid.org/#1}{\mbox{\scalerel*{
\begin{tikzpicture}[yscale=-1,transform shape]
\pic{orcidlogo};
\end{tikzpicture}
}{|}}}}

%% file: CustomDefs.tex
\def\im{\mathrm{Im}\, }
\def\re{\mathrm{Re}\, }
\def\diam{\mathrm{diam}\, }
\def\Hess{\mathrm{Hess}\, }
\def\rank{\mathrm{rank }\,  }
\def\rg{\mathrm{rg }\,  }
\def\d{\mathrm{d}}
\def\p{\partial }
\def\pr{\mathrm{pr}}
\def\D{\mathrm{D}}
\def\id{\mathrm{id}}
\def\Id{\mathrm{Id}}
\def\e{\epsilon}
\def\C{\mathbb{C}}
\def\Z{\mathbb{Z}}
\def\Q{\mathbb{Q}}
\def\N{\mathbb{N}}
\def\R{\mathbb{R}}
\def\K{\mathbb{K}}
\def\T{\mathbb{T}}
\def\diag{\mathrm{diag}\, }

%% file: sec_geodesics.tex
\section{Geodesics on submanifolds as optimal control problems}\label{sec:GeoOC}

To prepare the computation of cut loci, we formulate the geodesic equation on Riemannian submanifolds as variational problems and OCPs.

\subsection{Variational formulation}
Consider a submanifold $M \subset \R^n$ given as the zero level set $M=g(0)$ of a smooth function $g \colon \R^n \to \R^m$, where $m \le n$ and the Jacobian matrix of $g$ at each point in $M$ is of maximal rank. Let $\| \cdot \|$ denote the Euclidean norm in $\R^n$.
A curve $q \in \mathcal{C}^\infty([0,1],M) \subset \mathcal{C}^\infty([0,1],\R^n)$ is a geodesic on $M$ that connects $q(0)=q_0$, $q(1)=q_N$ for $q_0, q_N \in M$ if the length functional
\[
S(q) = \frac{1}{2}\int_0^1\|\dot{q}(t)\|^2 \d t
\]
is stationary at $q$ among all curves of $\mathcal{C}^\infty([0,1],M)$ connecting $q_0$ and $q_N$. More precisely,
\[
\delta S(q)(v) = \lim_{\e \to 0} \frac{1}{\e}(S(q+ \e v)-S(q)) =0
\]
for all $v \in \mathcal{C}^\infty([0,1],M)$ with $v(0)=0=v(1)$.
Equivalently, $q \in \mathcal{C}^\infty([0,1],\R^n)$ with $q(0)=q_0$, $q(1)=q_N$ is a geodesic on $M$ if there exists a Lagrangian multiplier $\lambda \in \mathcal{C}^\infty([0,1],\R^m)$ such that $(q,\lambda)$ is a stationary point of the extended functional
\[
\bar{S}(q,\lambda) = \int_0^1 \left(\frac{1}{2} \|\dot{q}(t)\|^2 - g(q(t))^\top \lambda(t) \right) \d t,
\]
i.e.\ $\delta \overline{S}(q,\lambda)(v,w) =0$ for all variations $v \in \{ v \in \mathcal{C}^\infty([0,1],\R^n)\, | \, v(0)=0=v(1) \}$ and $w \in \mathcal{C}^\infty([0,1],\R^m)$.
Here, $g(q(t))^\top$ denotes the transposition of $g(q(t))$.

Using partial integration and the fundamental theorem of variational calculus on the condition $\delta \overline{S}(q,\lambda)(v,w) =0$ shows that state and Lagrangian multiplier $(q,\lambda)$ constitute a stationary point of $\overline{S}$ if and only if the boundary conditions and constrained Euler-Lagrange equations
\begin{equation}\label{eq:ConstEL}
	\begin{split}
	\ddot{q} &= - g'(q)^\top \lambda,\;
	g(q)=0,\;
	q(0)=q_0,\, q(1)=q_N
	\end{split}
\end{equation}
are fulfilled.

\subsection{Optimal control formulation}

The variational principle $\delta \overline{S}=0$ or the constrained equations \eqref{eq:ConstEL} could be discretised directly. However, we would like to analyse which discretisation strategies for OCPs preserve qualitative aspects of the solution space to OCPs, which admit several extremal solutions. To generate an example which is simple on the one hand but rich enough on the other hand to exhibit the considered phenomena, we will cast the problem of finding geodesics as an OCP as follows.

%
%
A cost functional
\[
J(u,\lambda) = \int_0^1 \left( \frac{1}{2} \|u(t)\|^2 - g(q(t))^\top \lambda(t) \right) \d t
\]
is extremised among all controls $(u,\lambda) \in \mathcal{C}^\infty([0,1],\R^n\times \R^m)$ with associated states $q\in \mathcal{C}^\infty([0,1],\R^n)$ that fulfil the state equation
\[\dot{q} = u, \quad \text{subject to}\quad q(0)=q_0, q(1)=q_N.\]

\begin{prop}
To any optimal control $(u,\lambda) \in \mathcal{C}^\infty([0,1],$ $\R^n\times \R^m)$ and its associated state trajectory $q$ there exists a costate trajectory $p$ such that
\begin{equation}
	\label{eq:PMPODE}
	\begin{split}
		\dot{q}=u,\;
		\dot{p}=-g'(q)^\top \lambda
	\end{split}
\end{equation}
subject to the algebraic constraints
\begin{equation}
	\label{eq:PMPALG}
	g(q)=0, \quad p= u
\end{equation}
and the boundary conditions $q(0)=q_0$, $q(1)=q_N$.
\end{prop}

\begin{pf}
	By Pontryagin's principle \citep{Liberzon2012} there exists a scalar $p_0 \le 0$ and a costate trajectory $p \colon [0,1] \to \R^n$ with $(p_0,p) \not = (0,0)$ such that the optimal control $(u,\lambda)$ and its associated state trajectory $q$ fulfil Hamilton's equations
	\begin{equation}\label{eq:HamEQ}
		\dot{q} = \frac{\p H}{\p p}, \quad \dot{p} = -\frac{\p H}{\p q}
	\end{equation} 
	for the Hamiltonian
	\[H(q,p,u,\lambda)=p^\top u + p_0\left(\frac{1}{2} \|u\|^2 - g(q)^\top \lambda\right)\]
	and, as the optimal control variables $(u,\lambda)$ are assumed to take values in the open set $\R^n \times \R^m$, the optimality condition $\frac{\p H}{\p (u,\lambda)} = 0$ holds.	
	In particular
	\[
	0=\frac{\p H}{\p u} = p+p_0 u.
	\]
	If the abnormal multiplier $p_0$ is zero, then $p \equiv 0$ which contradicts the non-triviality condition $(p_0,p) \not = (0,0)$. After rescaling, if necessary, we can assume $p_0 =-1$. Now \eqref{eq:PMPODE} and \eqref{eq:PMPALG} are obtained from \eqref{eq:HamEQ} and the optimality condition $\frac{\p H}{\p (u,\lambda)} = 0$.
	 \hspace{\fill} \qed
\end{pf}

Relations \eqref{eq:PMPODE} and \eqref{eq:PMPALG} constitute first order necessary conditions for optimal controls $(u,\lambda) \in \mathcal{C}^\infty([0,1],\R^n\times \R^m)$. This yields the following constrained boundary value problem for optimal state and costate trajectories: 
\begin{equation}
	\label{eq:PMPbvp}
	\begin{split}
		\dot{q}&=p,\;
		\dot{p}=-g'(q)^\top \lambda\\
		g(q)&=0, \;
		q(0)=q_0,\, q(1)=q_N.
	\end{split}
\end{equation}
Notice that \eqref{eq:PMPbvp} recovers \eqref{eq:ConstEL}.

\section{Discretisation}\label{sec:Disc}

We now apply several popular discretisation schemes to the different formulations of the geodesic problem.

\subsection{Discrete Euler-Lagrange equations}

We discretise the action $\overline{S}$ to
\begin{equation}\label{eq:SDelta}
\overline{S}_\Delta = \frac{1}{2}\frac{\|q_{1}-q_0\|^2}{{\Delta t}}+\sum_{k=1}^{N-1} \left(\frac{1}{2}\frac{\|q_{k+1}-q_k\|^2}{{\Delta t}} - \Delta t g(q_k)^\top \lambda_{k}\right),
\end{equation}
where $\Delta t>0$ is a discretisation parameter. The values $(q_k,\lambda_k) \in \R^n \times \R^m$ are interpreted as an approximation to $(q(k \Delta t),\lambda(k \Delta t))$ for $k=1,\ldots,N-1$.
The values $q_0$ and $q_N$ are determined by the boundary condition.
An extremum $\{(q_k,\lambda_k)\}_{k=1}^{N-1}$ fulfils
\begin{equation*}
	\label{eq:DELComp}
	\begin{split}
		0&= \frac{\p \overline{S}_\Delta}{\p q_k} =	-\frac{q_{k+1}-2q_k + q_{k-1}}{\Delta t} - \Delta t g'(q_{k})^\top \lambda_k \\
		0&= \frac{\p \overline{S}_\Delta}{\p \lambda_k} = -g(q_{k})
	\end{split}
\end{equation*}
with $k=1,\ldots,N-1$. This induces a scheme in which $q_{k+1}$ can be computed from $q_k$ and $q_{k-1}$ by solving the $n+m$-dimensional system
\begin{equation}
	\label{eq:DEL}
	\begin{split}
0&= \frac{q_{k+1}-2q_k + q_{k-1}}{\Delta t} + \Delta t g'(q_{k})^\top \lambda_k, \;
0= g(q_{k+1})
	\end{split}
\end{equation}
for $q_{k+1}$ and $\lambda_k$. If $q_0$ and $q_N$ are given, then a collection of the formulas \eqref{eq:DEL} with $k=1,\ldots,N-1$ together with $g(q_1)=0$ constitutes an $(N-1)(n+m)$-dimensional system of nonlinear equations\footnote{The condition $g(q_{k+1})=0$ is removed from the last instance $k=N$ of \eqref{eq:DEL} as it is fulfilled by assumption.} which can be solved numerically with an iterative method. However, to reduce dimensionality {\em shooting methods} are usually preferred: 
for fixed $q_0$, an iteration of \eqref{eq:DEL} for $k=1,\ldots,N-2$ yields a map $q_1 \mapsto (q_{N-2},q_{N-1})$. A composition with the map $(q_{N-2},q_{N-1},\lambda_{N-1}) \mapsto q_N$ with
\[
q_N = 2q_{N-1}-q_{N-2} - \Delta t^2 g'(q_{N-1})^\top \lambda_{N-1}
\]
yields a map $\psi \colon (q_1,\lambda_{N-1}) \mapsto q_N$. To given $q_N$, the value $q_1$ can be computed numerically from the $n+m$-dimensional system $\psi(q_1,\lambda_{N-1})-q_N=0$, $g(q_1)=0$. Finally, obtain $\{q_k\}_{k=1}^{N-1}$ from \eqref{eq:DEL}. 

If, on the other hand, instead of the boundary values $q_0$ and $q_N$ an initial state $q_0$ and momentum $p_0$ (corresponding to tangential velocity) are given, then $q_1$ can be obtained by a discrete Legendre transformation \citep{MarsdenWestVariationalIntegrators} by solving
\[
p_0 = \frac{q_1-q_0}{\Delta t} + g'(q_0)^\top \tilde{\lambda}_0, \quad g(q_1)=0
\] 
for $q_1$ and the Lagrange multiplier $\tilde{\lambda}_0$. Afterwards, $\{q_k\}_{k=1}^{N}$ are obtained from \eqref{eq:DEL}. 



\subsection{Indirect method with symplectic discretisation}
We discretise \eqref{eq:PMPbvp} by a symplectic integrator, for instance, the symplectic Euler-Method:
\begin{align*}
	q_{k+1} &= q_k + \Delta t p_k\\
	p_{k+1} &= p_k -\Delta t g'(q_{k+1})^\top \lambda_{k+1},\;
	0 = g(q_{k+1})
\end{align*}

In the above scheme, the variables $p_k$ can be eliminated
such that
\[
\frac{q_{k+1}-2q_{k} + q_{k-1}}{\Delta t} +\Delta t g'(q_{k})^\top \lambda_{k} =0.
\]
The scheme is, therefore, equivalent to \eqref{eq:DEL}. The values for $p_k$ can be computed in a post-processing step, if required. 

\subsection{Indirect method with non-symplectic discretisation}
We discretise \eqref{eq:PMPbvp} by the (non-symplectic) explicit midpoint rule:
\begin{equation}\label{eq:RK2}
	\begin{split}
		q_{k+1} &= q_k + \Delta t \left(p_k - \frac{\Delta t}{2}g'(q_k)^\top \lambda_k\right), \; 0 = g(q_{k+1}) \\
		p_{k+1} &= p_k -\Delta t g'\left(q_{k}+\frac{\Delta t}{2}p_k\right)^\top \lambda_{k}
	\end{split}
\end{equation}

We will see later that this scheme is {\em not} equivalent to \eqref{eq:DEL} because it has different preservation properties.

\subsection{Karush–Kuhn–Tucker condition (direct method)}
When using the Karush–Kuhn–Tucker condition (KKT) for discretisation, we do {\em not} discretise the first order optimality conditions \eqref{eq:PMPbvp} but discretise the cost functional $J$ first and then derive first order optimality conditions for the discrete OCP (direct method).
For this, we first apply a numerical method to the state equation $\dot{q}=u$. The Euler method yields $q_{k+1}=q_k + \Delta t u_k$.
Next, the cost functional $J$ is discretised to
\[
J_\Delta = \sum_{k=0}^{N-1} \frac{1}{2} \|u_k\|^2 - g(q_k)^\top \lambda_k + \mu_k^\top(q_{k+1}-q_k-\Delta t u_k ),
\]
where the discretised state equation has been incorporated into the discrete action using Lagrange multipliers $\mu_k$. Then $J_\Delta$ is extremised. We obtain
\begin{align*}
	0=\frac{\p J_\Delta }{ \p \mu_k} &= q_{k+1}-q_k-\Delta t u_k, \; 
	0=\frac{\p J_\Delta }{ \p u_k} = u_k - \Delta t \mu_k\\
	0=\frac{\p J_\Delta }{ \p q_k} &= -g'(q_k)^\top \lambda_k - \mu_k + \mu_{k-1}, \;
	0=\frac{\p J_\Delta }{ \p \lambda_k} = -g(q_k)^\top. 
\end{align*}
for $k=1,\ldots,N-1$. Eliminating $u_k = \Delta t \mu_k$ and $\mu_k = (q_{k+1}-q_k)/\Delta t^2$ we can recover the scheme \eqref{eq:DEL}.

%% file: sec_experiment.tex
\section{Numerical Experiment}\label{sec:Ex}

\begin{figure}
	
\begin{center}
	\includegraphics[width=0.6\linewidth]{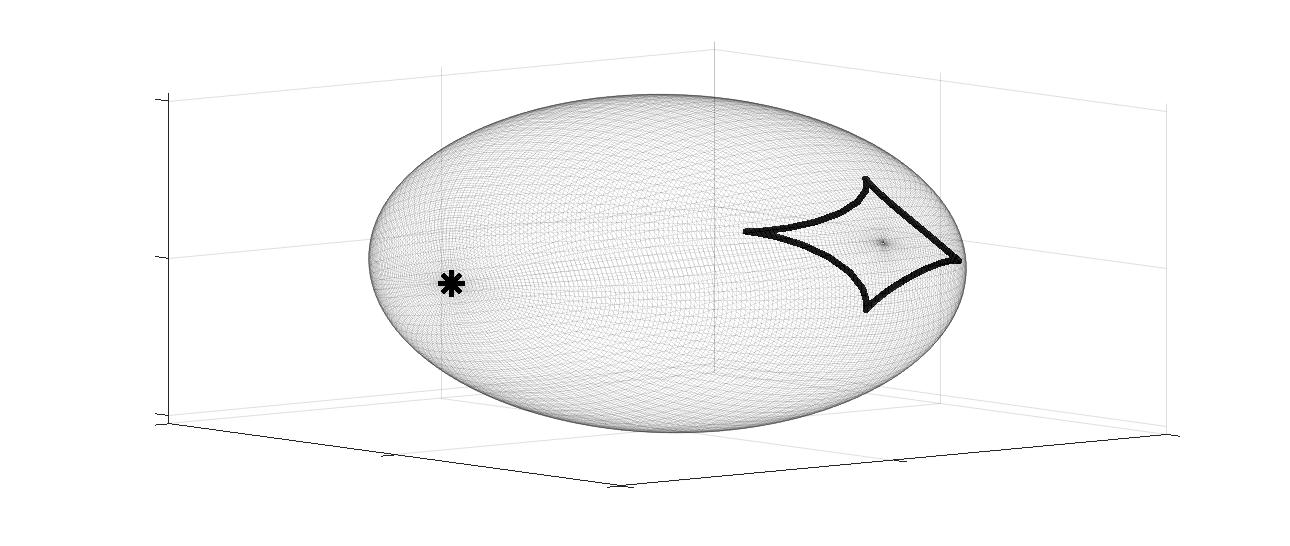}
\end{center}
	\caption{
		The cut locus $C$ on a 2-dimensional ellipsoid $E$ with respect to a  point $q_0 \in E$ marked by $\ast$. Denote the connected components of $E\setminus C$ by $E_0$ and $E_1$ such that $q_0 \in E_0$. 
		To any point $q_1 \in E_1$ there are three geodesics connecting $q_0$ with $q_1$ with length at most by $\frac{3}{2} \pi$, while there is only one such geodesic if $q_1 \in E_0$. 
		As $q_1$ is moved from $E_1$ to $E_0$ through a line of $C$, two of the geodesics merge and annihilate (fold bifurcation). If $q_1$ is moved through one of the four corners of $C$, three of the geodesics merge into one in a cusp bifurcation. If $q_0$ is in general position, $C$ always consists of four cusps connected by fold lines \citep{Itoh2004}.
	}\label{fig:2DLocus}
	
\end{figure}

The shortest path that connects two points on a (complete and connected) Riemannian manifold is always a geodesic. 
Moreover, each geodesic $\gamma\colon [0,T] \to M$ is locally length minimising, i.e.\ there exists an $\e>0$ such that for all $t \in [0,\e)$ the geodesic $\gamma$ is the shortest path connecting $\gamma(0)$ with $\gamma(t)$. If $\e$ is maximal with the property that $\gamma|_{[0,\e)}$ is length minimising, then $\gamma(\e)$ is a {\em conjugate point} to $\gamma(0)$. The set of all conjugate points to a reference point $q_0$ is referred to as the {\em cut locus to $q_0$} (see \citep{doCarmo} for exact definitions). The cut locus of a 2-dimensional ellipsoid in $\R^3$ is displayed in figure \ref{fig:2DLocus}.

A computation of a cut locus to a point $q_0$ on a Riemannian manifold $M$ given as a level set of a smooth function $g$ proceeds as follows.

\begin{itemize}
	\item
	Given a step-size $\Delta t>0$ and $N \in \N$ such that $1 = \Delta t N$ the scheme \eqref{eq:DEL} provides a map $M \ni q_1 \mapsto q_N \in M$.
	An open neighbourhood of $q_0 \in M$ can be identified with an open neighbourhood $\mathcal{O}$ of $0 \in \R^{m-n}$ such that the scheme induces a computable map $\phi \colon \mathcal O \to M$. The map $\phi$ corresponds to an expression of the classical geodesic exponential map in a chart.
	
	\item
	We compute the set of critical points $C_0 \subset \mathcal{O}$ of $\phi$. For this, we consider the determinant $\det \phi'$ of the Jacobian $\phi'$ of $\phi$, whereas the Jacobian of $\phi$ is computed using automatic differentiation. The critical set $C_0$ can be computed by evaluating $\phi$ on a mesh on $\mathcal{O}$ and then using a level-set method (such as {\tt contour} or {\tt isosurface} from Python's plotly package or MATLABs {\tt contourf} or {\tt isosurface}). Alternatively, $C_0$ can be computed by a level set continuation method such as pseudo-arclength continuation, if $C_0$ is one-dimensional, or by manifold continuation methods \citep{KrauskopfNumericalContinuation}.
	
	\item
	The critical set $C_0$ is mapped with $\phi$ to the set of critical values $C$, which is the cut locus.
	
	\item
	Additionally, highly degenerate points within $C_0$ (such as cusps, swallowtail points, etc.) can be computed using techniques such as those developed by one of the authors in \citep{PDEBifur}.
	
\end{itemize}

\begin{figure}
	\begin{center}
	\includegraphics[width=.45\linewidth]{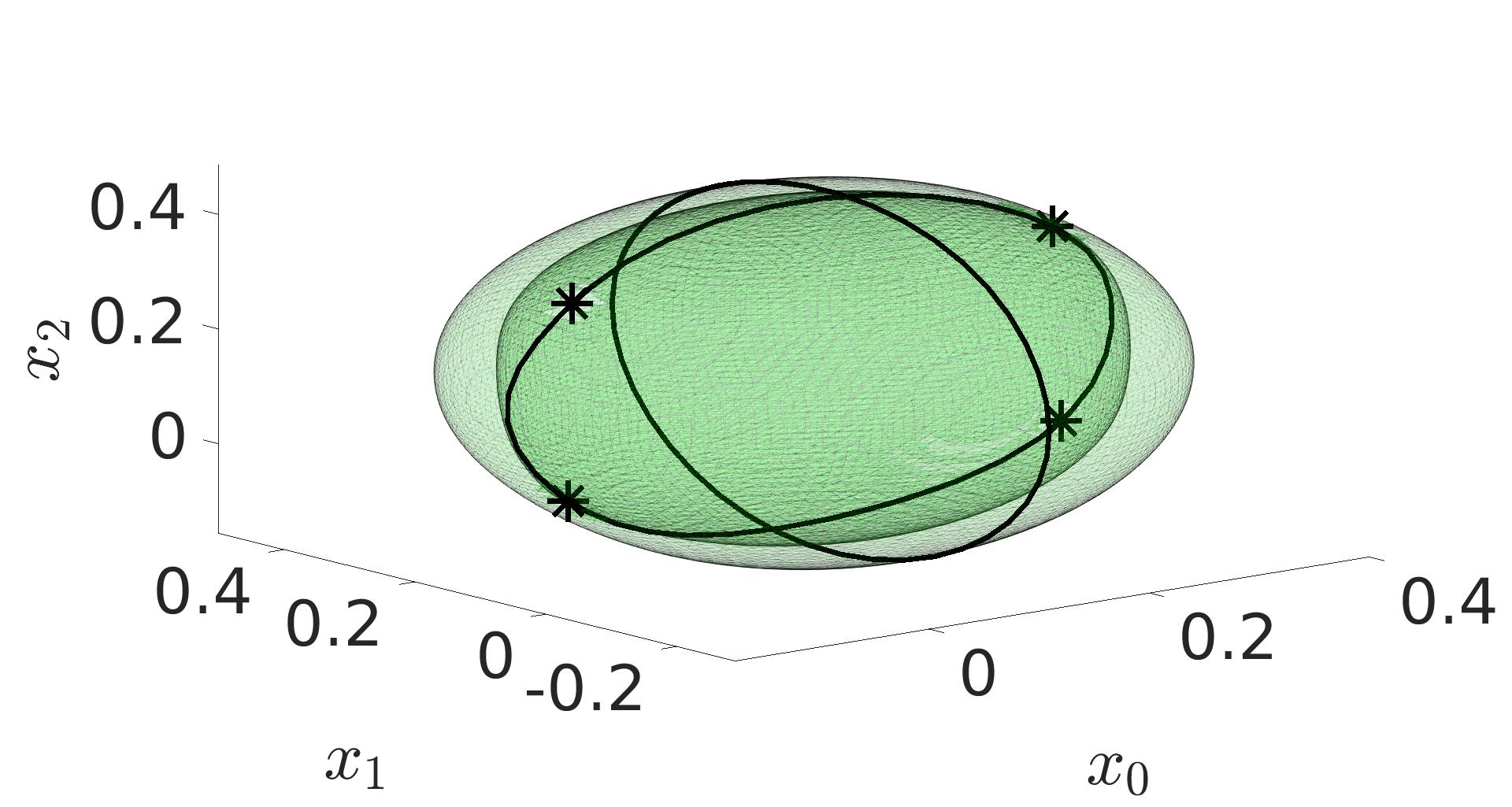}
	\includegraphics[width=.45\linewidth]{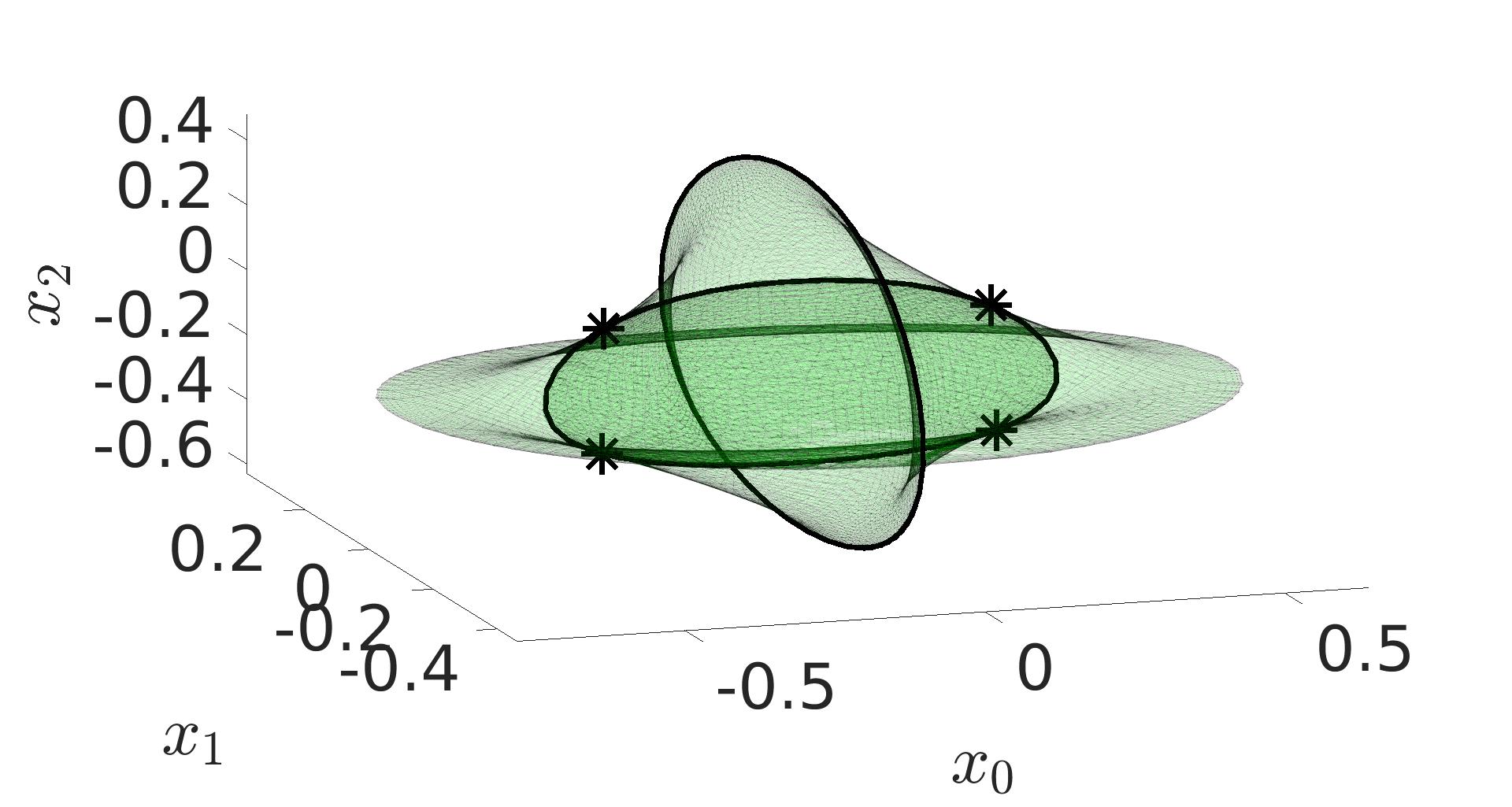}
	\end{center}
	\caption{The first figure shows the critical set $C_0$ and the second figure the conjugate locus of a 3-dimensional ellipsoid. Each point on the sheets corresponds to a (preimage of a) fold singularity. The solid lines correspond to (preimages of) lines of cusp singularities. At the points marked by $\ast$ there are umbilic singularities which we investigate closer in figure \ref{fig:Umbilic}. The cusp lines and locations of umbilic singularities have been computed using techniques from \citep{PDEBifur,obstructionPaper}}\label{fig:Locus3D}
\end{figure}

Figure \ref{fig:Locus3D} shows the conjugate locus with respect to a typical point of a three-dimensional ellipsoid considered as a submanifold of $\R^4$. Only the first three components $x_0,x_1,x_2$ are plotted. The last component $x_3$ can be recovered from $x_0,x_1,x_2$. A theoretical description of the structure of loci of high-dimensional ellipsoids has recently been given in \citep{Itoh2020}. Numerical computations on a normal form of a high-dimensional ellipsoid can be found in \citep{Joets1999}.

\begin{figure}
	\begin{center}
	\includegraphics[width=0.6\linewidth]{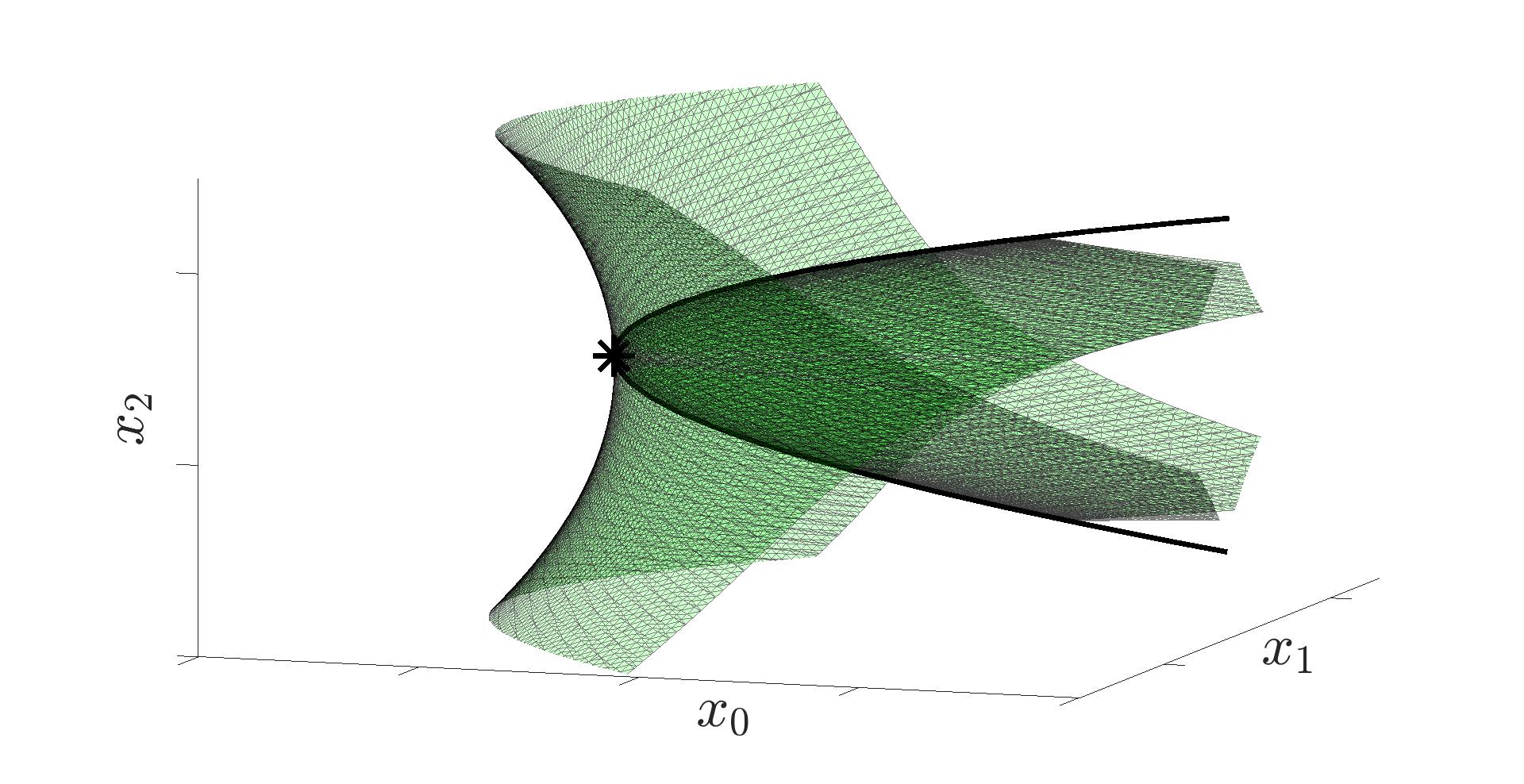}
	\includegraphics[width=0.45\linewidth]{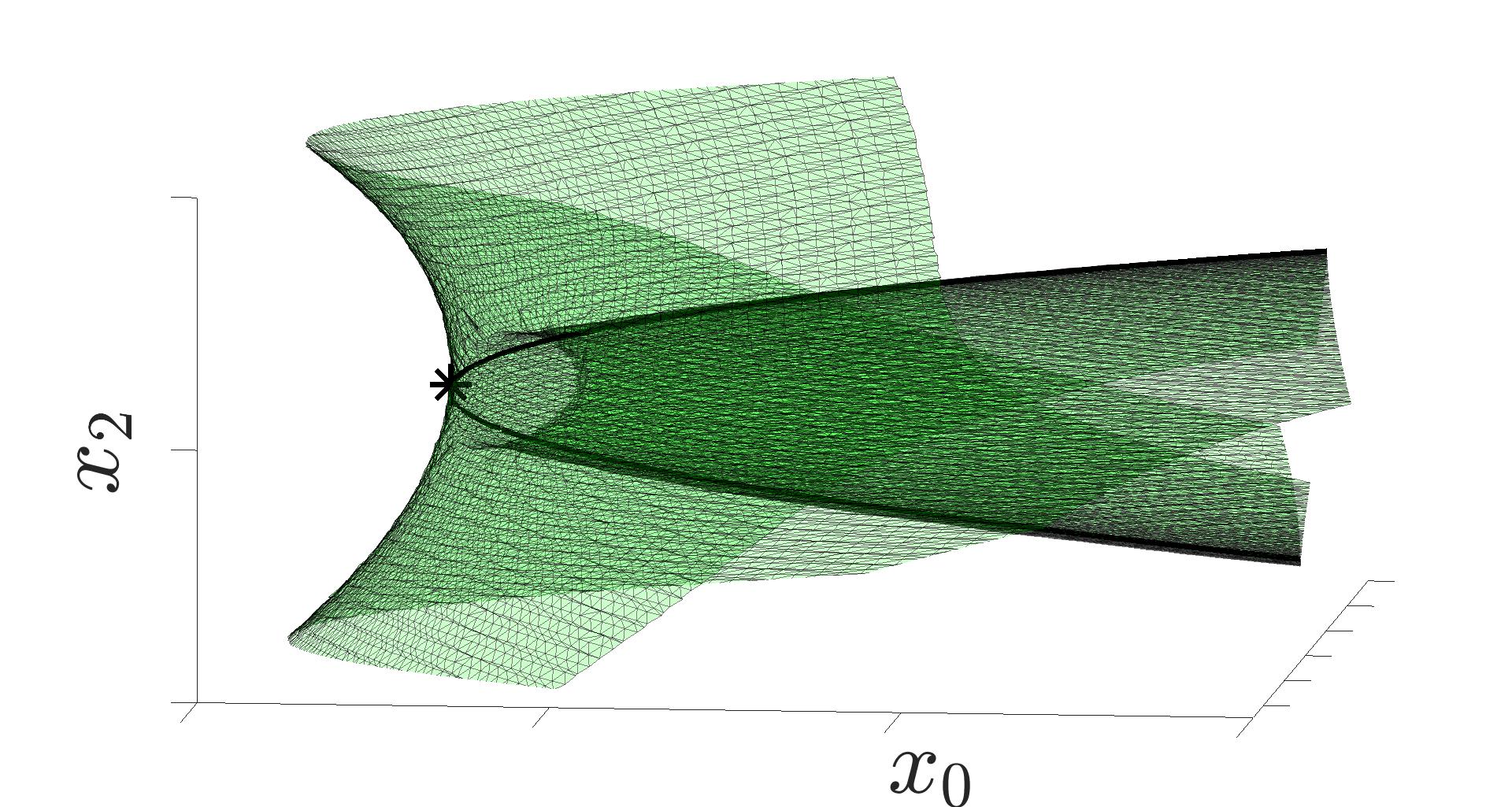}
	\includegraphics[width=0.45\linewidth]{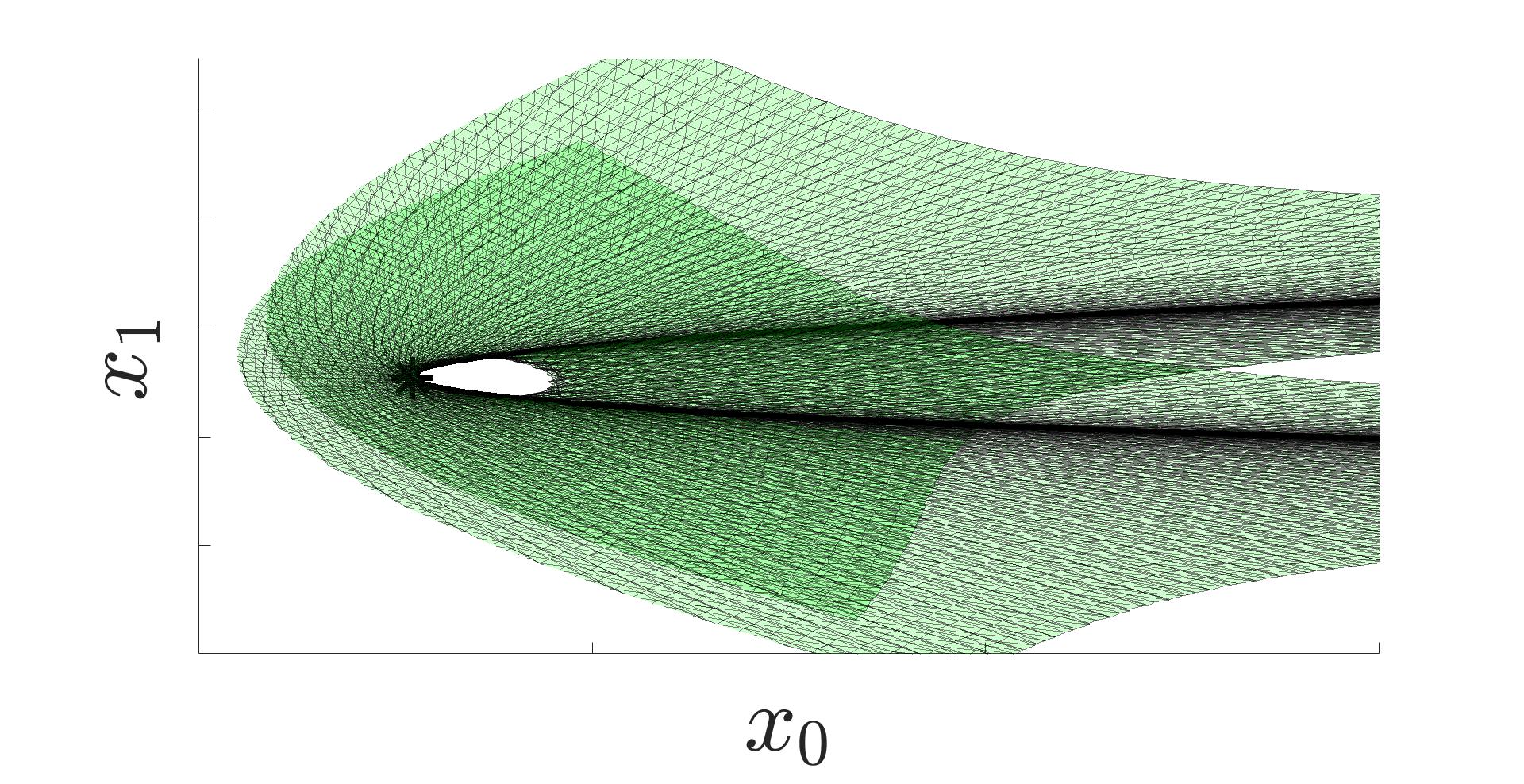}
	\caption{The plots show the conjugate locus close to an umbilic bifurcation point. While the top figure was obtained by the variational scheme \eqref{eq:DEL}, the other two plots correspond to an experiment with \eqref{eq:RK2}. 
		While the locus to \eqref{eq:DEL} correctly shows a hyperbolic umbilic bifurcation, the locus to \eqref{eq:RK2} is broken and contains an artificial hole and does {\em not} contain a hyperbolic umbilic singularity. Indeed, the most singular point marked by $\ast$ corresponds to a simple fold singularity rather than an umbilic singularity. See \cite{YoutubeSingularitiesAnimations} for animations that have been obtained by investigating normal forms of classical catastrophes.}\label{fig:Umbilic}
		\end{center}
\end{figure}

\begin{figure}
	\begin{center}
	\includegraphics[width=0.45\linewidth]{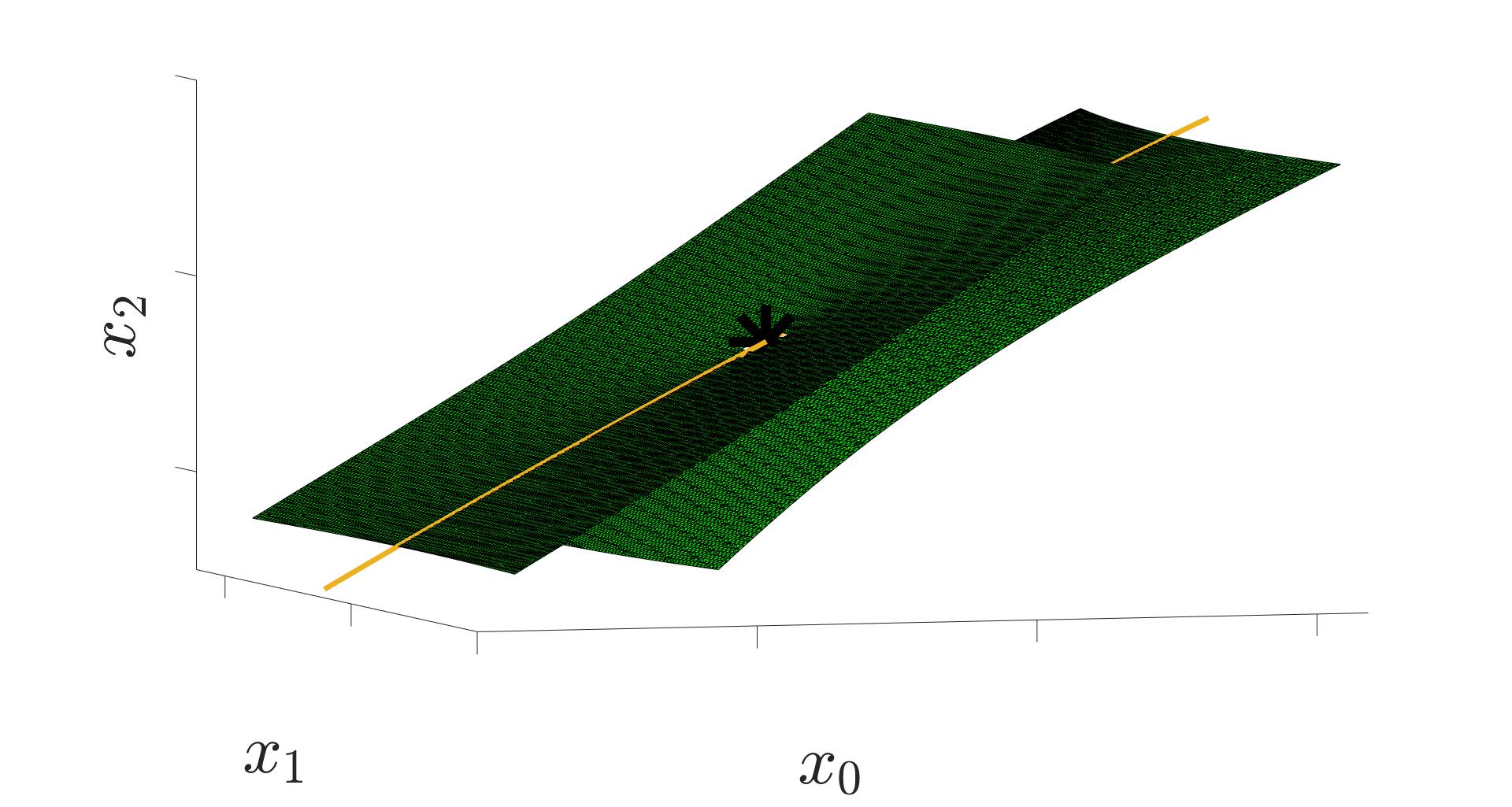}
	\includegraphics[width=0.45\linewidth]{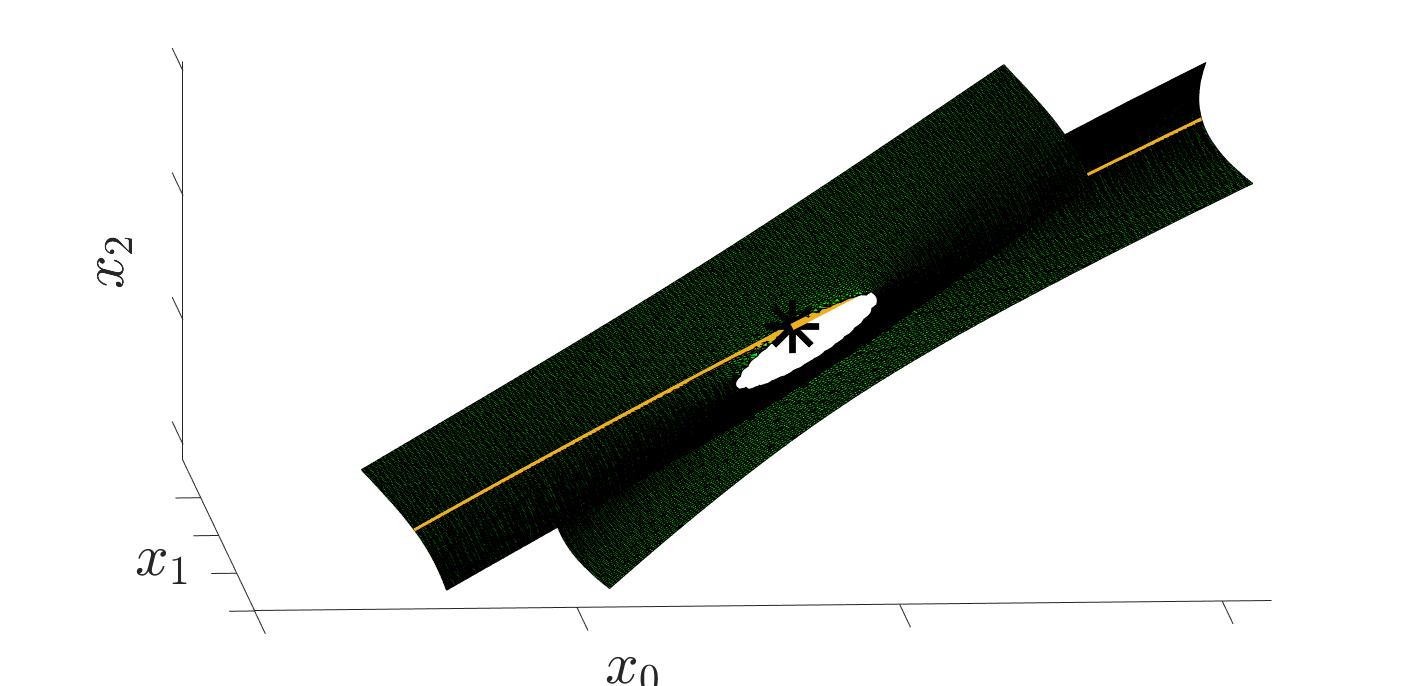}
	\caption{The plots show the preimage $C_0$ of the conjugate locus close to an umbilic bifurcation point. The left figure was obtained by the variational scheme \eqref{eq:DEL} and correctly shows two tangentially intersecting sheets. The right figure was obtained by the non-variational scheme \eqref{eq:RK2} and does not capture the situation well. }\label{fig:UmbilicPre}
		\end{center}
\end{figure}

Figure \ref{fig:Umbilic} shows (rotated) close-ups of the cut locus near one of the hyperbolic umbilic points.
Only the variational scheme \eqref{eq:DEL} correctly captures the hyperbolic umbilic bifurcation. The non-variational method \eqref{eq:RK2} breaks the bifurcation, contains an artificial hole, and two spurious highly degenenerate points at the locations where the line of cusp bifurcations touches the sheet of folds.
The structural error of the non-variational scheme can also be seen in the computed critical sets $C_0$, i.e.\ the preimage on the locus. While the first plot of figure \ref{fig:UmbilicPre} corresponding to the variational scheme correctly shows two sheets intersecting in one point, in the second plot two sheets connect along a circle. 
%
%
%
%
This demonstrates the significance of structure preservation for the computation of solutions to variational problems such as OCPs when bifurcations occur. Source code is available in \citep{LocusSoftware}.

%% file: sec_general.tex
\section{Theoretical considerations}\label{sec:Theo}

We will now develop a theoretical framework to predict which local bifurcations occur generically in parameter dependent OCPs and to explain the different behaviour of symplectic and non-symplectic discretisation schemes. Using Pontryagin's principle, the local bifurcation behaviour of regular solutions to necessary conditions of OCPs will be translated to the bifurcation behaviour of solutions to Hamiltonian boundary value problems, which are related to catastrophy theory \citep{bifurHampaper,PhDThesis}.
In a neighbourhood of a smooth solution the infinite-dimensional setting of OCPs is, thus, reduced to a setting to which finite-dimensional theory applies. 

Let the state space be given by a smooth manifold $M$ without boundary and let the control space $U$ and parameter space $\Lambda$ be open subsets $U \subset \R^l$, $\Lambda \subset \R^k$. For a given parameter $\mu \in \Lambda$ we consider the extremisation of the cost function $S$ from \eqref{eq:OCPGeneralIntro} 
among smooth controls $u \in \mathcal{C}^\infty([t_0,t_N],U)$ subject to the state equation \eqref{eq:stateEQ}.
The parameter $\mu$ is fixed during the extremisation process.
Here $q_0(\mu),q_N(\mu) \in M$ and for each $u \in U$ the function $f(\cdot,u;\mu) \colon M \to TM$ is a smooth vectorfield on $M$. Moreover, $L$ is smooth and all data depends smoothly on the parameter $\mu$. Let $\pi \colon T^\ast M \to M$ denote the cotangent bundle projection and let $\langle \cdot, \cdot \rangle$ be the pairing of $T^\ast M$ and $TM$. 
Assume that for $\mu^\ast \in \Lambda$ there exists a control $u^\ast \in \mathcal{C}^\infty([t_0,t_N],U)$ with associated state space trajectory $q^\ast \in \mathcal{C}^\infty([t_0,t_N],M)$ which extremises $S(\cdot;\mu^\ast)$. By Pontryagin's principle \citep{Barbero2008} there exists an abnormal multiplier $p_0^\ast \le 0$ and a costate trajectory $\alpha^\ast \colon [t_0,t_N] \to T^\ast M$ with $\pi(\alpha^\ast)=q^\ast$ and $(p_0^\ast,\alpha^\ast) \not = (0,0_{T^\ast M})$, where $0_{T^\ast M}$ is the zero section in $T^\ast M$, such that $\alpha^\ast$ is a trajectory of the Hamiltonian vectorfield $X_{H_{u^\ast,\mu^\ast}}$ defined by the Hamiltonian
$H \colon T^\ast M \times U \times \Lambda \to \R$,
$H (\alpha,u;\mu) = \langle \alpha, f(q,u;\mu) \rangle + p_0 L(q,u;\mu)$. 
In other words, the Hamiltonian boundary value problem
\begin{equation}\label{eq:HamBVP}
	\begin{split}
&\frac{\d}{\d t} \alpha(t) = X_{H_{u,\mu}}(\alpha(t))\\ &\pi(\alpha(t_0)) = q_0(\mu),\; \pi(\alpha(t_N)) = q_N(\mu)
	\end{split}
\end{equation}
is fulfilled for $(\alpha,u;\mu) = (\alpha^\ast,u^\ast;\mu^\ast)$. Here $X_{H_{u,\mu}}$ denotes the Hamiltonian vectorfield to the Hamiltonian $H_{u,\mu}(\alpha) = H(\alpha,u;\mu)$.
Moreover, since the controls $u$ take values in the open set $U \subset \R^l$, the optimality condition
\begin{equation}\label{eq:HamOpt}
	\frac{\p H}{\p u}(\alpha(t),u(t);\mu) = 0, \quad \forall t \in [t_0,t_N]
\end{equation}
holds true for $(\alpha,u;\mu) = (\alpha^\ast,u^\ast;\mu^\ast)$. 
The Hamiltonian boundary value problem \eqref{eq:HamBVP} and the optimality condition \eqref{eq:HamOpt} constitute first order necessary conditions for the OCP \eqref{eq:OCPGeneralIntro}.

\begin{thm}\label{thm:BifurCat}
	Assume that the Hessian matrix $\frac{\p H}{\p u^i \p u^j}$ is invertible. All stable, local bifurcations of smooth solutions to the first order necessary conditions with an abnormal multiplier $p_0 \not =0$ of the OCP \eqref{eq:OCPGeneralIntro} are governed by catastrophe theory.
	Direct discretisation methods for OCPs as well as symplectic discretisation methods for the first order necessary conditions preserve all stable bifurcations.
\end{thm}

\begin{rem}
	
	\begin{itemize}		
		
		
		\item
		As the constraint considered in the OCP for the computation of a conjugate locus of a 3-dimensional ellipsoid $M$ is holonomic, the problem is equivalent to a problem covered by theorem \ref{thm:BifurCat}, where the ellipsoid is the manifold $M$. The parameters enters only in the boundary condition $q(t_N) = q_N(\mu)$, where $q_N \colon \Lambda \to M$ is a local chart of $M$ with $\Lambda \subset \R^3$.

		\item
		Stability in this context means that the bifurcations are persistent under small perturbations of the necessary condition \eqref{eq:HamBVP} within the class of Hamiltonian boundary value problems, or, slightly more generally, the class of boundary value problems for symplectic maps. Admissible perturbations of the OCP \eqref{eq:OCPGeneralIntro} include small perturbations of the state equation, the Lagrangian $L$ as well as the boundary condition for $q$ (whereas care needs to be taken if the perturbed boundary condition is allowed to involve $\dot{q}$, which leads to the notion of Lagrangian boundary conditions \citep{PhDThesis}).
		
		\item
		The non-degeneracy assumption on $\frac{\p H}{\p u^i \p u^j}$ only needs to hold in a tubular neighbourhood of the graph of a solution $(\alpha^\ast, u^\ast;\mu^\ast)$ for the theorem to hold close to $(q^\ast, u^\ast; \mu^\ast)$.
		
		\item
		If the state space $M$ or the control space $U$ contain boundaries, then the theorem can be applied locally by restricting to a tubular neighbourhood of the graph of $(q^\ast, u^\ast;\mu^\ast)$ if the image of the curves $q^\ast$, $u^\ast$ do not intersect with the boundaries.

	\end{itemize}
	
\end{rem}

\begin{pf}
	Let $(\alpha^\ast, u^\ast; \mu^\ast)$ be a solution to the first order necessary conditions \eqref{eq:HamBVP} and \eqref{eq:HamOpt}.
	Applying the implicit function theorem to $\frac{\p H}{\p u}(\alpha,u;\mu) =0$, there exists a unique function $\bar{u}$ depending on $(t,\alpha;\mu)$ such that $\bar{u}(t,\alpha^\ast(t);\mu^\ast)=u^\ast(t)$ and $\frac{\p H}{\p u}(\alpha,u(t,\alpha;\mu);\mu) =0$. The function $\bar{u}$ is defined on the Cartesian product of a tubular neighbourhood of the graph of $(\alpha^\ast, u^\ast)$ and an open neighbourhood of $\mu^\ast$ in $\Lambda$.
	Close to $(\alpha^\ast, u^\ast; \mu^\ast)$ the first order necessary conditions \eqref{eq:HamBVP} and \eqref{eq:HamOpt} are, therefore, equivalent to the (time-dependent) Hamiltonian boundary value problem
		$\frac{\d}{\d t} \alpha(t) = X_{H_{\bar{u},\mu}}(\alpha(t))$,
		$\pi(\alpha(t_0)) = q_0(\mu),\; \pi(\alpha(t_N)) = q_N(\mu)$
	The flow map of $X_{H_{\bar{u},\mu}}$ is a symplectic map such that we obtain a family of boundary value problems for symplectic maps $\{\phi_\mu\}_{\mu \in \Lambda^\ast}$, where $\Lambda^\ast \subset \Lambda$ is an open neighbourhood of $\mu^\ast$.
	The local bifurcation behaviour of solutions to such systems was related to catastrophe theory in \citep{bifurHampaper,PhDThesis}. 
	In an indirect discretisation method of \eqref{eq:OCPGeneralIntro}, Hamilton's equations \eqref{eq:HamBVP} are discretised and the control is (as before) obtained from the optimality condition \eqref{eq:HamOpt}. If a symplectic integrator is used in \eqref{eq:HamBVP}, then we obtain a family of boundary value problems for symplectic maps $\bar{\phi}_\mu$ close to $\phi_\mu$. A stable, catastrophe type bifurcation of the boundary value problem for $\phi_\mu$ is present in the nearby problem for $\bar{\phi}$ as well for sufficiently small discretisation parameters by the definition of stability.
	As direct discretisation methods for OCPs correspond to indirect methods with symplectic integration \citep{ober2011DMOC}, the conclusion also holds for direct methods.	
	 \hspace{\fill} \qed
\end{pf}

As shown in \citep{numericalPaper}, $D$-series bifurcations, such as hyperbolic umbilic bifurcations, are stable bifurcations in families of boundary value problems for symplectic maps but unstable in more general classes of boundary value problems. Other bifurcations, such as fold, cusp, which belong to the $A$-series, are also stable in wider classes of boundary value problems. Therefore, the correspondence of regular OCP and symplectic boundary value problems, provided by the proof of theorem \ref{thm:BifurCat}, explains our observations from the numerical example that hyperbolic umbilic bifurcations are preserved when using a structure preserving discretisation schemes, while the sheets of fold bifurcations and lines of cusp singularities persist even if variational or symplectic structure is destroyed under discretisation.

\section{Conclusions and discussion}\label{sec:Conclusion}

Solutions to first order necessary conditions of families of optimal control problems can undergo bifurcations as parameters are varied. Under regularity assumptions we showed that local bifurcations which are persistent under small perturbations of the family of OCPs are exactly the classical catastrophes. Moreover, to preserve all stable bifurcations under discretisation, either direct discretisation methods for OCPs or indirect methods in combination with symplectic integrators can be used. If, however, discretisation methods are used which are not structure preserving, then certain bifurcations, such as $D$-series bifurcations, break. A preservation of bifurcations is necessary when computing bifurcation diagrams to determine in which parameter ranges how many first order optimal solutions exist. 



%% file: ifacconf.bbl
\begin{thebibliography}{20}
\providecommand{\natexlab}[1]{#1}
\providecommand{\url}[1]{\texttt{#1}}
\providecommand{\urlprefix}{URL }
\expandafter\ifx\csname urlstyle\endcsname\relax
  \providecommand{\doi}[1]{doi:\discretionary{}{}{}#1}\else
  \providecommand{\doi}{doi:\discretionary{}{}{}\begingroup
  \urlstyle{rm}\Url}\fi

\bibitem[{Arnold et~al.(1998)Arnold, Goryunov, Lyashko, and Vasil'ev}]{Arnold1}
Arnold, V.I., Goryunov, V.V., Lyashko, O.V., and Vasil'ev, V.A. (1998).
\newblock \emph{Critical Points of Functions}, 10--50.
\newblock Springer Berlin Heidelberg, Berlin, Heidelberg.
\newblock \doi{10.1007/978-3-642-58009-3_1}.

\bibitem[{Barbero-Li{\~{n}}{\'{a}}n and
  Mu{\~{n}}oz-Lecanda(2008)}]{Barbero2008}
Barbero-Li{\~{n}}{\'{a}}n, M. and Mu{\~{n}}oz-Lecanda, M. (2008).
\newblock Geometric approach to {P}ontryagin's maximum principle.
\newblock \emph{Acta Applicandae Mathematicae}, 108(2), 429--485.
\newblock \doi{10.1007/s10440-008-9320-5}.
\newblock \urlprefix\url{https://doi.org/10.1007%2Fs10440-008-9320-5}.

\bibitem[{Chyba et~al.(2009)Chyba, Hairer, and
  Vilmart}]{SymplecticityOptimalControl}
Chyba, M., Hairer, E., and Vilmart, G. (2009).
\newblock The role of symplectic integrators in optimal control.
\newblock \emph{Optimal Control Applications and Methods}, 30(4), 367--382.
\newblock \doi{10.1002/oca.855}.
\newblock
  \urlprefix\url{https://onlinelibrary.wiley.com/doi/abs/10.1002/oca.855}.

\bibitem[{Flaherty and do~Carmo(1992)}]{doCarmo}
Flaherty, F. and do~Carmo, M. (1992).
\newblock \emph{Riemannian Geometry}.
\newblock Mathematics: Theory \& Applications. Birkh{\"a}user Boston.
\newblock \doi{10.1007/978-1-4757-2201-7}.
\newblock \urlprefix\url{https://doi.org/10.1007%2F978-1-4757-2201-7}.

\bibitem[{Itoh and Kiyohara(2004)}]{Itoh2004}
Itoh, J. and Kiyohara, K. (2004).
\newblock The cut loci and the conjugate loci on ellipsoids.
\newblock \emph{manuscripta mathematica}, 114(2), 247--264.
\newblock \doi{10.1007/s00229-004-0455-z}.
\newblock \urlprefix\url{https://doi.org/10.1007/s00229-004-0455-z}.

\bibitem[{Itoh and Kiyohara(2020)}]{Itoh2020}
Itoh, J. and Kiyohara, K. (2020).
\newblock The structure of the conjugate locus of a general point on ellipsoids
  and certain liouville manifolds.
\newblock \emph{Arnold Mathematical Journal}.
\newblock \doi{10.1007/s40598-020-00153-9}.
\newblock \urlprefix\url{https://doi.org/10.1007%2Fs40598-020-00153-9}.

\bibitem[{Joets and Ribotta(1999)}]{Joets1999}
Joets, A. and Ribotta, R. (1999).
\newblock {Caustique de la surface ellipsoïdale à trois dimensions}.
\newblock \emph{Experimental Mathematics}, 8(1), 49 -- 55.
\newblock \doi{em/1047477111}.
\newblock \urlprefix\url{https://doi.org/em/1047477111}.

\bibitem[{Kogan(1986)}]{Kogan1986}
Kogan, J. (1986).
\newblock \emph{Bifurcation of Extremals in Optimal Control}.
\newblock Springer Berlin Heidelberg.
\newblock \doi{10.1007/bfb0077060}.
\newblock \urlprefix\url{https://doi.org/10.1007%2Fbfb0077060}.

\bibitem[{Krauskopf et~al.(2007)Krauskopf, Osinga, and
  Gal{\'{a}}n-Vioque}]{KrauskopfNumericalContinuation}
Krauskopf, B., Osinga, H.M., and Gal{\'{a}}n-Vioque, J. (eds.) (2007).
\newblock \emph{Numerical Continuation Methods for Dynamical Systems}.
\newblock Springer Netherlands.
\newblock \doi{10.1007/978-1-4020-6356-5}.
\newblock \urlprefix\url{https://doi.org/10.1007%2F978-1-4020-6356-5}.

\bibitem[{Kreusser et~al.(2020)Kreusser, McLachlan, and Offen}]{PDEBifur}
Kreusser, L.M., McLachlan, R.I., and Offen, C. (2020).
\newblock Detection of high codimensional bifurcations in variational {PDE}s.
\newblock \emph{Nonlinearity}, 33(5), 2335--2363.
\newblock \doi{10.1088/1361-6544/ab7293}.
\newblock \urlprefix\url{https://doi.org/10.1088/1361-6544/ab7293}.

\bibitem[{Liberzon(2012)}]{Liberzon2012}
Liberzon, D. (2012).
\newblock \emph{The Maximum Principle}, 102--155.
\newblock Princeton University Press.
\newblock \doi{10.2307/j.ctvcm4g0s}.
\newblock \urlprefix\url{http://www.jstor.org/stable/j.ctvcm4g0s}.

\bibitem[{Marsden and West(2001)}]{MarsdenWestVariationalIntegrators}
Marsden, J.E. and West, M. (2001).
\newblock Discrete mechanics and variational integrators.
\newblock \emph{Acta Numerica}, 10, 357–514.
\newblock \doi{10.1017/S096249290100006X}.
\newblock \urlprefix\url{https://dx.doi.org/10.1017/S096249290100006X}.

\bibitem[{McLachlan and Offen(2018{\natexlab{a}})}]{bifurHampaper}
McLachlan, R.I. and Offen, C. (2018{\natexlab{a}}).
\newblock Bifurcation of solutions to {H}amiltonian boundary value problems.
\newblock \emph{Nonlinearity}, 31(6), 2895--2927.
\newblock \doi{10.1088/1361-6544/aab630}.
\newblock \urlprefix\url{https://doi.org/10.1088/1361-6544/aab630}.

\bibitem[{McLachlan and Offen(2018{\natexlab{b}})}]{obstructionPaper}
McLachlan, R.I. and Offen, C. (2018{\natexlab{b}}).
\newblock Hamiltonian boundary value problems, conformal symplectic symmetries,
  and conjugate loci.
\newblock \emph{New Zealand Journal of Mathematics (NZJM)}, 48, 83--99.
\newblock
  \urlprefix\url{http://nzjm.math.auckland.ac.nz/index.php/Hamiltonian_Boundary_Value_Problems\%2C_Conformal_Symplectic_Symmetries\%2C_and_Conjugate_Loci}.

\bibitem[{McLachlan and Offen(2019)}]{numericalPaperShort}
McLachlan, R.I. and Offen, C. (2019).
\newblock Symplectic integration of boundary value problems.
\newblock \emph{Numerical Algorithms}, 81(4), 1219--1233.
\newblock \doi{10.1007/s11075-018-0599-7}.
\newblock \urlprefix\url{https://doi.org/10.1007/s11075-018-0599-7}.

\bibitem[{McLachlan and Offen(2020)}]{numericalPaper}
McLachlan, R.I. and Offen, C. (2020).
\newblock Preservation of bifurcations of {H}amiltonian boundary value problems
  under discretisation.
\newblock \emph{Foundations of Computational Mathematics (FoCM)}, 20,
  1363--1400.
\newblock \doi{10.1007/s10208-020-09454-z}.
\newblock \urlprefix\url{https://doi.org/10.1007/s10208-020-09454-z}.

\bibitem[{Ober-Bl{\"o}baum et~al.(2011)Ober-Bl{\"o}baum, Junge, and
  Marsden}]{ober2011DMOC}
Ober-Bl{\"o}baum, S., Junge, O., and Marsden, J.E. (2011).
\newblock Discrete mechanics and optimal control: an analysis.
\newblock \emph{ESAIM: Control, Optimisation and Calculus of Variations},
  17(2), 322--352.
\newblock \doi{10.1051/cocv/2010012}.
\newblock \urlprefix\url{https://dx.doi.org/10.1051/cocv/2010012}.

\bibitem[{Offen(2019)}]{YoutubeSingularitiesAnimations}
Offen, C. (2019).
\newblock Singularities animations.
\newblock
  \url{https://www.youtube.com/playlist?list=PLIp-UrijLTJ5m-3ZASHPurIkehiBuW_sO}.
\newblock Accessed 2021-05-02.

\bibitem[{Offen(2020)}]{PhDThesis}
Offen, C. (2020).
\newblock \emph{Analysis of {H}amiltonian boundary value problems and
  symplectic integration ({D}octoral {T}hesis)}.
\newblock Massey University.
\newblock \doi{10.13140/RG.2.2.34063.61607}.
\newblock \urlprefix\url{http://dx.doi.org/10.13140/RG.2.2.34063.61607}.

\bibitem[{Offen(2021)}]{LocusSoftware}
Offen, C. (2021).
\newblock Release v1.0 of {GitHub} repository {Christian-Offen/ConjugateLocus}.
\newblock \doi{10.5281/zenodo.4562664}.
\newblock \urlprefix\url{https://doi.org/10.5281/zenodo.4562664}.

\end{thebibliography}
